\documentclass{article}

\title{\LARGE \textbf{Large cycles in 4-connected graphs}}
\author{M.Zh. Nikoghosyan and Zh.G. Nikoghosyan\footnote{G.G. Nicoghossian up to 1997}}

\begin{document}

\maketitle

\begin{abstract}
Every 4-connected graph $G$ with minimum degree $\delta$ and connectivity $\kappa$ either contains a cycle of length at least $4\delta-\kappa-4$ or every longest cycle in $G$ is a dominating cycle.\\

\end{abstract}

We consider only finite undirected graphs without loops or multiple edges. Let $n$ denote the order, $\delta$ the minimum degree, $\kappa$ the connectivity and $c$ the circumference (the length of a longest cycle) of a graph $G$. A cycle $C$ is a Hamilton cycle if $|C|=n$ and is a dominating cycle if every edge of $G$ has a vertex in common with $C$. A cycle $C$ is said to be a $D_3$-cycle if every path of length at least 2 has a vertex in common with $C$.

In 2008, Yamashita [3] obtained a degree sum condition for dominating cycles which yields the following.\\

\noindent\textbf{Theorem A} [3]. Let G be a 3-connected graph. If $\delta\geq(n+\kappa+3)/4$, then any longest cycle in $G$ is a dominating cycle.\\

In this paper we prove, in fact, the reverse version of Theorem A.\\

\noindent\textbf{Theorem 1}. Let $G$ be a 4-connected graph. Then either $c\geq4\delta-\kappa-4$ or every longest cycle in $G$ is a dominating cycle.\\

In order to prove Theorem 1, we need the following result due to Jung [2].\\

\noindent\textbf{Theorem B} [2]. Let $G$ be a 4-connected graph. Then either $c\geq4\delta-8$ or every longest cycle in $G$ is a $D_3$-cycle.\\

A good reference for any undefined terms is [1]. The set of vertices of a graph $G$ is denoted by $V(G)$ and the set of edges by $E(G)$. For $S$ a subset of $V(G)$, we denote by $G\backslash S$ the maximum subgraph of $G$ with vertex set $V(G)\backslash S$. For a subgraph $H$ of $G$ we use $G\backslash H$ short for $G\backslash V(H)$. We denote by $N(x)$ the neighborhood of a vertex x in a graph $G$ with $d(x)=|N(x)|$. Furthermore, for a subgraph $H$ of $G$ and $X\subseteq V(G)$, we define $N_H(X)=N(X)\cap V(H)$. 

Paths and cycles in a graph $G$ are considered as subgraphs of $G$. If $Q$ is a path or a cycle, then the length of $Q$, denoted by $|Q|$, is $|E(Q)|$. We write a cycle $C$ with a given orientation by $\overrightarrow{C}$. For $x,y\in V(C)$, we denote by $x\overrightarrow{C}y$, or sometimes by $C[x,y]$, the subpath of $C$ in the chosen direction from $x$ to $y$. For $C[x^+,y^+]$ we also write $C(x,y)$. For $x\in V(C)$, we denote the $h$-th successor and the $h$-th predecessor of $x$ on $\overrightarrow{C}$ by $x^{+h}$ and $x^{-h}$, respectively. We abbreviate $x^{+1}$ and $x^{-1}$by $x^+$ and $x^-$, respectively. For $X\subset V(C)$, we define $X^{+h}=\{x^{+h}|x\in X\}$ and $X^{-h}=\{x^{-h}|x\in X\}$. 

Henceforth, we use the following notation. Let $G$ be a 4-connected graph and $C$ be a $D_3$-cycle in $G$ with $x_1x_2\in E(G\backslash C)$. We denote
$$
R=N_C(x_1)\cup N_C(x_2), \quad   M=N_C(x_1)\cap N_C(x_2),
$$
$$
A=R\backslash M, \quad A_1=N_C(x_1)\backslash M, \quad A_2=N_C(x_2)\backslash M, 
$$
$$
Y=R\cup R^+\cup M^{+2}.
$$  

\noindent\textbf{Lemma 1}. $|Y|\geq 2d(x_1)+2d(x_2)-|M|-4\geq 4\delta-|M|-4.$\\

\noindent\textbf{Proof}. Since $C$ is extreme, $R,R^+$ and $M^{+2}$ are pairwise disjoint. Observing that $R=A\cup M$, we get
$$
|Y|=|R|+|R^+|+|M^{+2}|=2|R|+|M|=3|M|+2|A|.
$$

Next, since $|A_i|=d(x_i)-|M|-1$ and $|A|=|A_1|+|A_2|$,
$$
|Y|\geq 2d(x_1)+2d(x_2)-|M|-4\geq 4\delta-|M|-4.  \qquad \Delta
$$

\noindent\textbf{Lemma 2}. Let $P$ be a longest $(x,y)$-path, having only $y$ in common with $C$. If $y\in R^+$, then either $|P|=0$ or $c\geq 4\delta-\kappa-3$, and if $y\in M^{+2}$, then either $|P|\leq 1$ or $c\geq 4\delta-\kappa-3$.\\

\noindent\textbf{Proof}. Let $\xi_1,...,\xi_t$ be the elements of $R$, occuring on $\overrightarrow{C}$ in a consecutive order. Assume w.l.o.g. that $y\in \{\xi_1^+,\xi_1^{+2}\}$ and choose $w\in N_C(x)\backslash\{y\}$. Let $w\in V(\xi^+_i\overrightarrow{C}\xi_{i+1})$ for some $i\in \{1,...,t\}$. Let $Q$ be a longest path connecting $\xi_i$ to $\xi_1$ and passing through $\{x_1,x_2\}$. Put $C^{\prime}=\xi_iQ\xi_1\overleftarrow{C}wxPy\overrightarrow{C}\xi_i$.\\

\textbf{Case 1}. $y\in R^+$.

Assume that $|P|\geq1$. Since $C$ is a $D_3$-cycle, $1\leq|P|\leq2$.\\

\textbf{Case 1.1.} $P=xy$.

Since $|C^\prime|\leq|C|$, we have $|\xi_i\overrightarrow{C}w|\geq4$ if $\xi_i\in M$ and $|\xi_i\overrightarrow{C}w|\geq3$ if $\xi_i\in A$. It means that
$$
(N^-(x)\backslash\{y^-,y^+\})\cap Y=\emptyset.
$$ 

Observing also that $d(x)\geq|M|+1$ and using Lemma 1, we get 
$$
c\geq|Y|+|N^-(x)\backslash\{y^-,y^+\}|\geq4\delta+d(x)-|M|-6\geq4\delta-\kappa-3.
$$

\textbf{Case 1.2.} $P=xzy$ for some vertex $z$.

Observing that $(N^-_C(x)\backslash\{y^-\})\cap Y=\emptyset$ and $|N_C(x)|\geq|N(x)|-1$, we can argue as in Case 1.1.\\ 

\textbf{Case 2}. $y\in M^{+2}$.

Assume that $|P|\geq2$. Since $C$ is a $D_3$-cycle, we have $|P|=2$, i.e. $P=xzy$ for some $z\in V(G)$. Due to $|P|=2$, we can obtain $(N^-_C(x)\backslash \{y^-\})\cap Y=\emptyset$ and further we can argue as in Case 1.1. \qquad $\Delta$ \\

\noindent\textbf{Lemma 3}. Let $S$ be a minimum cut-set of $G$. Then  either $c\geq4\delta-\kappa-3$ or $\{x_1,x_2\}\cap S=\emptyset$ for each $x_1x_2\in E(G\backslash C)$.\\

\noindent\textbf{Proof}. Choose a longest cycle $C$ such that $|V(C)\cap S|$ is as great as possible and let $x_1x_2\in E(G\backslash C)$ with $\{x_1,x_2\}\cap S\neq \emptyset$. Let $\xi_1,...,\xi_t$ be the elements of $R$, occuring on $\overrightarrow{C}$ in a consecutive order. Since $C$ is extreme, $(R^+\cup M^{+2})\cap R=\emptyset$. Further, since $|V(C)\cap S|$ is maximum, $M_1^{+3}\cap R=\emptyset$. Observing also that $|M_2|\leq\kappa-1$ and using Lemma 1, we get
$$
c\geq |Y|+|M_1^{+3}|\geq4\delta-|M|-4+|M_1|=4\delta-|M_2|-4\geq4\delta-\kappa-3. \qquad \Delta
$$

\noindent\textbf{Proof of Theorem 1}. Let $G$ be a 4-connected graph, $S$ be a minimum cut-set in $G$ and $H_1,...,H_h$ be the components of $G\backslash S$. If $c\geq4\delta-8$, then we are done, since $4\delta-8\geq4\delta-\kappa-4$. Otherwise, by Theorem B, every longest cycle in $G$ is a $D_3$-cycle. Let $C$ be any longest cycle and $x_1x_2\in E(G\backslash C)$. Assume w.l.o.g. that $x_1x_2\in V(H_1)$. By Lemma 3, $\{x_1,x_2\}\cap S=\emptyset$. Abbreviate, $V_1=V(H_1)\cup S$. Assume first that $Y\subseteq V_1$. By Lemma 1,
$$
|V(C)\cap V_1|\geq|Y|\geq (2d(x_1)+d(x_2))+d(x_2)-|M|-4\geq3\delta-3.
$$

If $V(H_2)\subseteq V(C)$, then $|V(C\cap H_2)|\geq\delta-\kappa+1$ and
$$
c\geq|V(C\cap V_1)|+|V(C\cap H_2)|\geq 4\delta-\kappa-2.
$$

Otherwise, we choose $y\in V(H_2\backslash C)$. Since $|N_C(y)|\geq\delta-1$, we have $|V(C\cap H_2)|\geq |N_C(y)|-|S|$ and
$$
c\geq|V(C\cap V_1)|+|V(C\cap H_2)|\geq4\delta-\kappa-4.
$$

Now let $Y\not\subseteq V_1$. Since $R\subseteq V_1$, we have $R^+\cup M^{+2}\not\subseteq V_1$. \\

\textbf{Case 1}. $R^+\cap V(H_2)\neq\emptyset$.

Let $y\in R^+\cap V(H_2)$. By Lemma 2, $N(y)\subseteq V(C)$ and by standard arguments, $N(y)\cap (R^+\cap M^{+2})=\emptyset$. Since $|N(y)\cap Y|=|N(y)\cap R|\leq \kappa$, we have by Lemma 1,
$$
c\geq |Y|+|N(y)|-|N(y)\cap Y|\geq|Y|+\delta-\kappa
$$
$$
\geq (2d(x_1)+d(x_2)+\delta)+(d(x_2)-|M|)-\kappa-4\geq 4\delta-\kappa-3.
$$

\textbf{Case 2}. $R^+\cup V(H_2)=\emptyset$.

We have $R\cup R^+\subseteq V_1$ and $M^{+2}\cap V(H_2)\neq \emptyset$. Then it is easy to see that for each $y\in M^{+2}\cap V(H_2)$, $|N(y)\cap Y|\leq \kappa$. If $N(y)\subseteq V(C)$, then $c\geq |Y|+|N(y)\backslash Y|\geq |Y|+|N(y)|-\kappa$ and we can argue as in Case 1. Let $N(y)\not\subseteq V(C)$ and $z\in N(y)\backslash V(C)$. By Lemma 2, $N(z)\subseteq V(C)$. If $z\in V(H_2)$, then by standard arguments, $N(z)\cap (R^+\cup M^{+2})=\emptyset$. Hence 
$$
c\geq |Y|+|N(y)\backslash Y|=|Y|+|N(y)|-|N(y)\cap R|\geq |Y|+|N(y)|-\kappa
$$
and we can argue as in Case 1. Let $z\not\in V(H_2)$. Then we can assume that $N(y)\backslash V(C)\subseteq S$. Set $D=N(y)\cap (S\backslash V(C))$. Since $|N(y)\cap Y|\leq \kappa-|D|$ and $|N(y)\cap V(C)|\geq \delta-|D|$, we have
$$
c\geq |Y|+|(N(y)\cap V(C))\backslash Y|
$$
$$
=|Y|+|N(y)\cap V(C)|-|N(y)\cap Y|\geq |Y|+\delta-\kappa
$$
and again we can argue as in Case 1.     \qquad $\Delta$


\begin{thebibliography}{9}

\bibitem{}	J.A. Bondy and U.S.R. Murty, Graph Theory with Applications. Macmillan, London and Elsevier, New York (1976).

\bibitem{}	H.A. Jung, Long cycles in graphs with moderate connectivity, in: R.Bodendick, R.Henn (Eds.), Topics in Combinatorics and Graph Theory, Phisika Verlag, Heidelberg (1990) 765-778 .

\bibitem{}	T. Yamashita, Degree Sum and Connectivity Conditions for Dominating Cycles, Discrete Math. 308 (9) (2008) 1620-1627.



\end{thebibliography}
\end{document}